\documentclass[11pt]{amsart}

\usepackage{graphicx}
\usepackage[top=26mm, bottom=26mm, left=28mm, right=28mm]{geometry}
\usepackage{multirow}
\usepackage{bm}
\usepackage{color}

\newcommand{\mtx}[1]{\mathsf{#1}}

\def\({\left(}
\def\){\right)}
\def\[{\left[}
\def\]{\right]}

\def\xb{\boldsymbol{x}}

\def\nb{\boldsymbol{n}}

\newtheorem{remark}{Remark}

\newcommand{\vct}[1]{\bm{\mathsf{#1}}}
\newcommand{\pvct}[1]{\bm{#1}}
\newcommand{\pxx}{\pvct{x}}
\newcommand{\pyy}{\pvct{y}}

\begin{document}

\begin{center}
\textbf{An efficient and highly accurate solver for multi-body acoustic scattering problems involving rotationally symmetric scatterers}

\vspace{4mm}

\textit{\small S.~ Hao, P.G.~Martinsson,  P.~Young}

\vspace{4mm}

\begin{minipage}{0.9\textwidth}\small
\noindent\textbf{Abstract:}
A numerical method for solving the equations modeling acoustic scattering
is presented.
The method is capable of handling several dozen scatterers, each of which
is several wave-lengths long, on a personal work station. Even for geometries
involving cavities, solutions accurate to seven digits or better were obtained.
The method relies on a Boundary Integral Equation formulation of the scattering
problem, discretized using a high-order accurate Nystr\"om method. A hybrid
iterative/direct solver is used in which a local scattering matrix for each
body is computed, and then GMRES, accelerated by the Fast Multipole Method,
is used to handle reflections between the scatterers. The main limitation of
the method described is that it currently applies only to scattering bodies that
are rotationally symmetric.
\end{minipage}
\end{center}


\section{Introduction}
The manuscript presents a robust and highly accurate numerical method for
modeling frequency domain acoustic scattering on a domain external to a
group of scatterers in three dimensions. The solver is designed for the
special case where each scatterer is rotationally symmetric, and relies
on a Boundary Integral Equation (BIE) formulation of the scattering problem.


The contribution of the manuscript is to combine several recently
developed techniques to obtain a solver capable of solving scattering
problems on complex multibody geometries in three dimensions to seven
digits of accuracy or more. In particular, the solver is capable of
resolving domains involving cavities such as, e.g., the geometry shown
Figure \ref{fig:cavity_8}(a).

The solution technique proposed involves the following steps:

\begin{enumerate}
\item \textit{Reformulation.} The problem is written mathematically as a BIE on
the surface of the scattering bodies using the ``combined field''
formulation \cite{1998_colton_kress_inverse,2012_nedelec}. See Section \ref{sec:BIE} for details.

\item \textit{Discretization.} The BIE is discretized using the Nystr\"om method based on a high-order
accurate composite Gaussian quadrature rule. Despite the fact that the kernel
in the BIE is singular, high accuracy can be maintained using the correction
techniques of \cite{2001_rokhlin_kolm,2013_martinsson_nystrom_final}.
Following \cite{Rizzo:79a}, we exploit the rotational symmetry of each body to decouple
the local equations as a sequence of equations defined on a generating contour
\cite{journal/jasa/80/4/10.1121/1.393817,Soenarko:93a,Kuijpers:97a,Wang:97a,Tsinopoulos:99a}.
This dimension reduction technique requires an
efficient method for evaluating the fundamental solution of the Helmholtz equation
in cylindrical coordinates (the so called ``toroidal harmonics''); we use
the technique described in \cite{2012_martinsson_axisymmetric}.
See Section \ref{sec:discretization} for details.

\item \textit{Iterative solver.} The dense linear system resulting from the Nystr\"om discretization of the BIE
is solved using the iterative solver GMRES \cite{1986_saad_gmres}, combined with a block-diagonal
pre-conditioner, as in, e.g., \cite[Sec.~6.4]{2012_ho_greengard_fastdirect}. This pre-conditioner
exploits that a highly accurate discrete approximation to the scattering matrix for each individual
scatterer can be computed efficiently. See Section \ref{sec:iteration} for details.

\item \textit{Fast matrix-vector multiplication.}
The application of the coefficient matrix in the iterative solver is acclerated
using the Fast Multipole Method (FMM) \cite{rokhlin1987}, specifically the version
for the Helmholtz equation developed by Gimbutas and Greengard \cite{2011_gimbutas_greengard_FMMLIB3D}.

\item \textit{Skeletonization.} In situations where the individual scatterers are not
packed very tightly, the number of degrees of freedom in the global system can be
greatly reduced by exploiting rank deficiencies in the off-diagonal blocks of the
coefficient matrix. Specifically, we use a variation of the scheme introduced in
\cite{lowrank}, and further developed in \cite{Martinsson:04a}.
Randomized methods are used to accelerate the computation of low-rank approximations
to large dense matrices \cite{MR2806637}.
See Section \ref{sec:accelerate} for details.
\end{enumerate}

The present work draws on several recent papers describing techniques for multibody
scattering, including \cite{2012_ho_greengard_fastdirect}, which applies a very similar
technique to acoustic scattering in two dimensions. \cite{journals/jcphy/GimbutasG13}
addresses the harder problem of electro-magnetic scattering in 3D (as opposed to the
acoustic scattering considered here), but uses classical scattering matrices expressed
in spherical harmonics. This is a more restrictive frame-work than the one used in
\cite{2012_ho_greengard_fastdirect} for problems in 2D, and in the present work for
problems in 3D. The more general model for a compressed scattering matrix that we use
here allows for larger scatterers to be handled, and also permits it to handle
scatterers closely packed together. For a deeper discussion of different ways of
representing compressed scattering matrices, see \cite{bremer2013high}.

To describe the asymptotic cost of the method presented, let $m$ denote the number of
scatterers, let $n$ denote the total number of discretization nodes on a single
scatterer and let $I$ denote the number of iterations required in our pre-conditioned
iterative solver to achieve convergence. The cost of building all local scattering
matrices is then $O(mn^{2})$, and the cost of solving the linear system consists of
the time $T_{\rm FMM}$ required for applying the coefficient matrices using the FMM,
and the time $T_{\rm precond}$ required for applying the block-diagonal preconditioner.
These scale as $T_{\rm FMM} \sim Imn$ and $T_{\rm precond} \sim Imn^{3/2}$
(cf.~Remark \ref{remark:costmatvec}), but for
practical problem sizes, the execution time is completely dominated by the FMM.
For this reason, we implemented a ``skeletonization'' compression scheme
\cite{lowrank} that reduces the cost of executing the FMM from $Imn$ to $Imk$,
where $k$ is a numerically determined ``rank of interaction''. We provide numerical
examples in Section \ref{sec:num} that demonstrate that when the scatterers
are moderately well separated, $k$ can by smaller than $n$ by one or two orders of
magnitude, leading to dramatic practical acceleration.



\section{Mathematical formulation of the scattering problem}
\label{sec:BIE}

\begin{figure}
\setlength{\unitlength}{1mm}
\begin{picture}(160,55)
\put(10,15){\includegraphics[height=40mm]{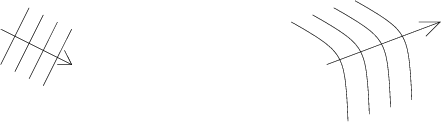}}
\put(00,21){\textit{incident field} $v$}
\put(128,14){\textit{scattered field} $u$}
\put(45,00){\includegraphics[height=45mm]{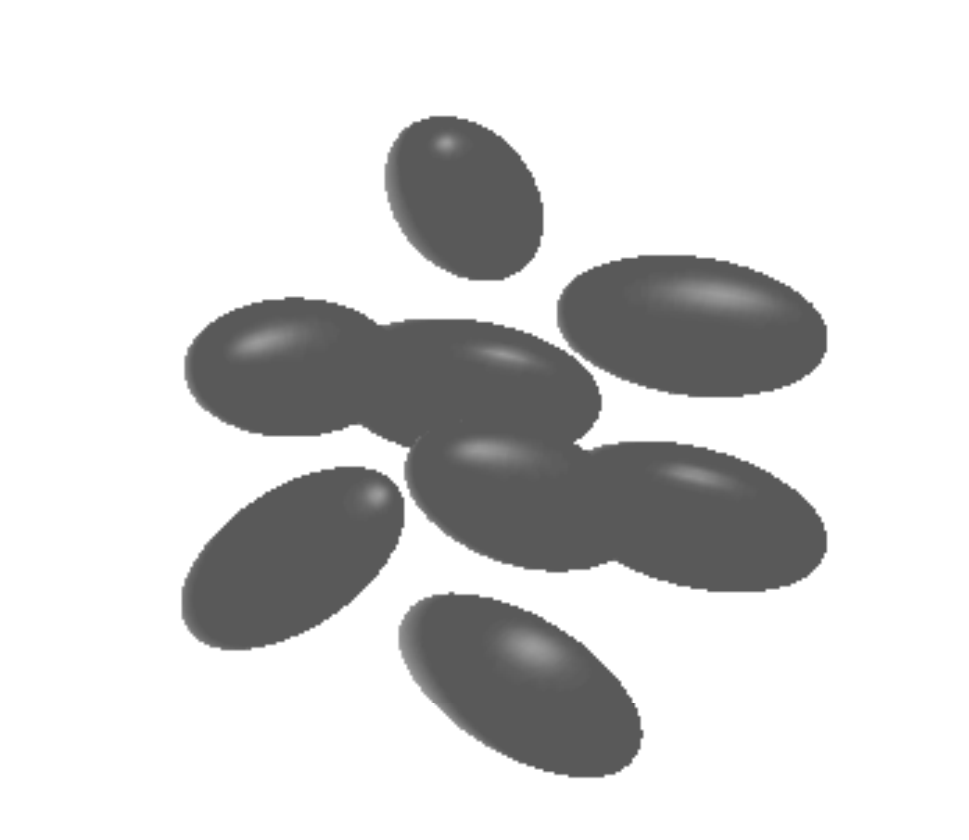}}
\end{picture}
\caption{Geometry of scattering problem. An incident field $v$ propagates in a medium with constant
wave-speed and hits a scattering surface $\Gamma = \bigcup_{p=1}^{m}\Gamma_{p}$ (shown for $m=8$).
A charge distribution $\sigma$ is induced on the surface $\Gamma$ and generates an outgoing field $u$.}
\label{fig:scattering_geometry}
\end{figure}

Let $\{\Gamma_{p}\}_{p=1}^{m}$
denote a collection of $m$ smooth, disjoint, rotationally symmetric surfaces
in $\mathbb{R}^{3}$, let  $\Gamma = \cup_{p=1}^{m} \Gamma_{p}$ denote
their union, and let $\Omega$ denote the domain exterior to $\Gamma$.
Our task is to compute the ``scattered field'' $u$ generated by an incident field
$v$ that hits the scattering surface $\Gamma$, see Figure \ref{fig:scattering_geometry}.
For concreteness, we consider the so called ``sound-soft'' scattering problem
\begin{equation}
\label{eq:basic}
\left\{\begin{aligned}
-\Delta u(\pxx) - \kappa^{2}u(\pxx) =&\ 0 &&\pxx \in \Omega^{\rm c},\\
u(\pxx) =&\ -v(\pxx) && \pxx \in \Gamma,\\
\frac{\partial u(\pxx)}{\partial r} - i\kappa u(\pxx) =&\ O(1/r) &&r:=|\pxx| \rightarrow\infty.
\end{aligned}\right.
\end{equation}
We assume that the ``wave number'' $\kappa$ is a real non-negative number. It
is known \cite{1998_colton_kress_inverse} that (\ref{eq:basic}) has a unique solution
for every incoming field $v$.

Following standard practice, we reformulate (\ref{eq:basic}) as second kind Fredholm
Boundary Integral Equation (BIE) using a so called ``combined field technique'' \cite{1998_colton_kress_inverse,2012_nedelec}.
We then look for a solution $u$ of the form
\begin{equation}
\label{eq:intHelm1}
u(\pxx) =
\int_{\Gamma}
G_{\kappa}(\pxx,\pxx')\,\sigma(\pxx')\,dA(\pxx'), \hspace{1em} \pxx \in \Omega^{c},
\end{equation}
where $G_{\kappa}$ is a combination of the single and double layer kernels,
\begin{equation}
\label{eq:Gkappa}
G_{\kappa}(\pxx,\pxx') =
\frac{\partial \phi_{\kappa}(\pxx,\pxx')}{\partial \nb(\pxx')} + i\kappa\,\phi_{\kappa}(\pxx,\pxx')
\end{equation}
and where $\phi_{\kappa}$ is the free space fundamental solution
\begin{equation}
\phi_{\kappa}(\pxx,\pxx') =\frac{e^{i \kappa |\pxx - \pxx'|}}{4\pi |\pxx - \pxx'|}.
\end{equation}
Equation (\ref{eq:intHelm1}) introduces a new unknown function $\sigma$, which we refer to as a
``boundary charge distribution''. To obtain an equation for $\sigma$, we take the limit
in (\ref{eq:intHelm1}) as $\pxx$ approaches the boundary $\Gamma$, and find that $\sigma$ must
satisfy the integral equation
\begin{equation}
\label{eq:intHelm2}
\frac{1}{2}\sigma(\pxx) +
\int_{\Gamma}G_{\kappa}(\pxx,\pxx')\,\sigma(\pxx')\,dA(\pxx') = -v(\pxx), \hspace{1em} \pxx \in \Gamma.
\end{equation}
The combined field equation (\ref{eq:intHelm2}) is known to be a second kind Fredholm equation
whenever $\Gamma$ is smooth. Like the orignal boundary value problem (\ref{eq:basic}), it is
known to be well posed for every $\kappa$, see \cite[Theorem.~3.9]{1998_colton_kress_inverse}, \cite[Sec.~3.2.2]{2012_nedelec} (in particular,
it does not suffer from the problem of ``artificial resonances'' that plague many alternative formulations).

\section{Discretization of rotationally symmetric scattering bodies}
\label{sec:discretization}

In Section \ref{sec:BIE} we formulated the scattering problem as the BIE (\ref{eq:intHelm2})
defined on the scattering surface $\Gamma$. In this section, we show how to discretize
(\ref{eq:intHelm2}) to obtain a system of linear algebraic equations $\mtx{A}\vct{\sigma}=-\vct{v}$.
We use a Nystr\"om technique that combines high accuracy, and (relative) ease of implementation.
Section \ref{sec:nystrom} gives a general overview of the Nystr\"om method,
Section \ref{sec:nystromsingle} describes how rotational symmetry can be exploited to relatively
easily discretize a single body to high order, and then Section \ref{sec:nystrommulti} describes
how to generalize the procedure to a multibody scattering problem.

\subsection{Nystr\"om discretization}
\label{sec:nystrom}
The Nystr\"om method provides a way of discretizing a BIE on a surface $\Gamma$ from a
quadrature rule for the surface that is valid for smooth functions.
To illustrate, suppose that we are given \textit{nodes} $\{\pxx_{i}\}_{i=1}^{n}$ and
\textit{weights} $\{w_{i}\}_{i=1}^{n}$ such that
\begin{equation}
\label{eq:basicquad}
\int_{\Gamma}\varphi(\pxx)\,dS(\pxx) \approx
\sum_{i=1}^{n} \varphi(\pxx_{i})\,w_{i},
\qquad\mbox{for}\
\varphi\
\mbox{smooth}.
\end{equation}
The idea is then to first use the discretization nodes $\{\pxx_{i}\}_{i=1}^{n}$ as collocation points;
in other words, we require that
\begin{equation}
\label{eq:collocated}
\frac{1}{2}\sigma(\pxx_{i}) +
\int_{\Gamma}G_{\kappa}(\pxx_{i},\pxx')\,\sigma(\pxx')\,dA(\pxx') = -v(\pxx_{i}), \hspace{1em} i = 1,2,3,\dots,n.
\end{equation}
Next, suppose that we can somehow (this can require some work) construct an $n\times n$ matrix $\mtx{A}$ such that
for any sufficiently smooth function $\varphi$, the integral in (\ref{eq:collocated}) can be approximated from
the function values $\{\sigma(\pxx_{i})\}_{i=1}^{n}$
\begin{equation}
\label{eq:nystromquad}
\frac{1}{2}\sigma(\pxx_{i}) +
\int_{\Gamma}G_{\kappa}(\pxx_{i},\pxx')\,\sigma(\pxx')\,dA(\pxx')
\approx
\sum_{j=1}^{n}\mtx{A}(i,j)\,\sigma(\pxx_{j})
\qquad\mbox{for}\ \sigma\ \mbox{smooth}.
\end{equation}
Then a system of $n$ equations for the $n$ unknowns $\{\sigma(\pxx_{i})\}_{i=1}^{n}$ is obtained by inserting
the approximation (\ref{eq:nystromquad}) into (\ref{eq:collocated}). Specifically, given a data vector $\vct{v} \in \mathbb{C}^{n}$
given by $\vct{v}(i) = v(\pxx_{i})$, we seek to determine a vector $\vct{\sigma} \in \mathbb{C}^{n}$ of approximations
$\vct{\sigma}(i) \approx \sigma(\pxx_{i})$ by solving the linear system
\begin{equation}
\label{eq:nystrom}
\sum_{j=1}^{n}\mtx{A}(i,j)\,\vct{\sigma}(j) =
-\vct{v}(i), \hspace{1em} i = 1,2,3,\dots,n.
\end{equation}

The task of constructing a matrix $\mtx{A}$ such that (\ref{eq:nystromquad}) holds
is complicated by the fact that the kernel $G_{\kappa}(\pxx,\pxx')$ has a singularity
as $\pxx' \rightarrow \pxx$. Had this not been the case, one could simply have applied
the rule (\ref{eq:basicquad}) to the integral in (\ref{eq:collocated}) to obtain
\begin{equation}
\label{eq:nystromnaive}
\mtx{A}(i,j) = G_{\kappa}(\pxx_{i},\pxx_{j})\,w_{j}.
\end{equation}
In Sections \ref{sec:nystromsingle} and \ref{sec:nystrommulti} we will describe how to
construct a basic quadrature rule $\{\pxx_{i},w_{i}\}_{i=1}^{n}$ that is suitable for the
geometry under consideration, and also how to construct a matrix $\mtx{A}$ such that
(\ref{eq:nystromquad}) holds to high accuracy despite the singular kernel. It turns out to be possible to do so
while having almost all elements of $\mtx{A}$ given by the simple formula
(\ref{eq:nystromnaive}) --- only matrix elements $\mtx{A}(i,j)$ for which
$||\pxx_{i}-\pxx_{j}||$ is ``small'' need to be modified. As we will see in
Section \ref{sec:iteration}, this will greatly help when forming fast algorithms for
evaluating the matrix-vector product $\vct{\sigma} \mapsto \mtx{A}\vct{\sigma}$.


\subsection{A single rotationally symmetric scatterer}
\label{sec:nystromsingle}
\begin{figure}[htbp]
	\centering
	\includegraphics[width = 0.3\linewidth]{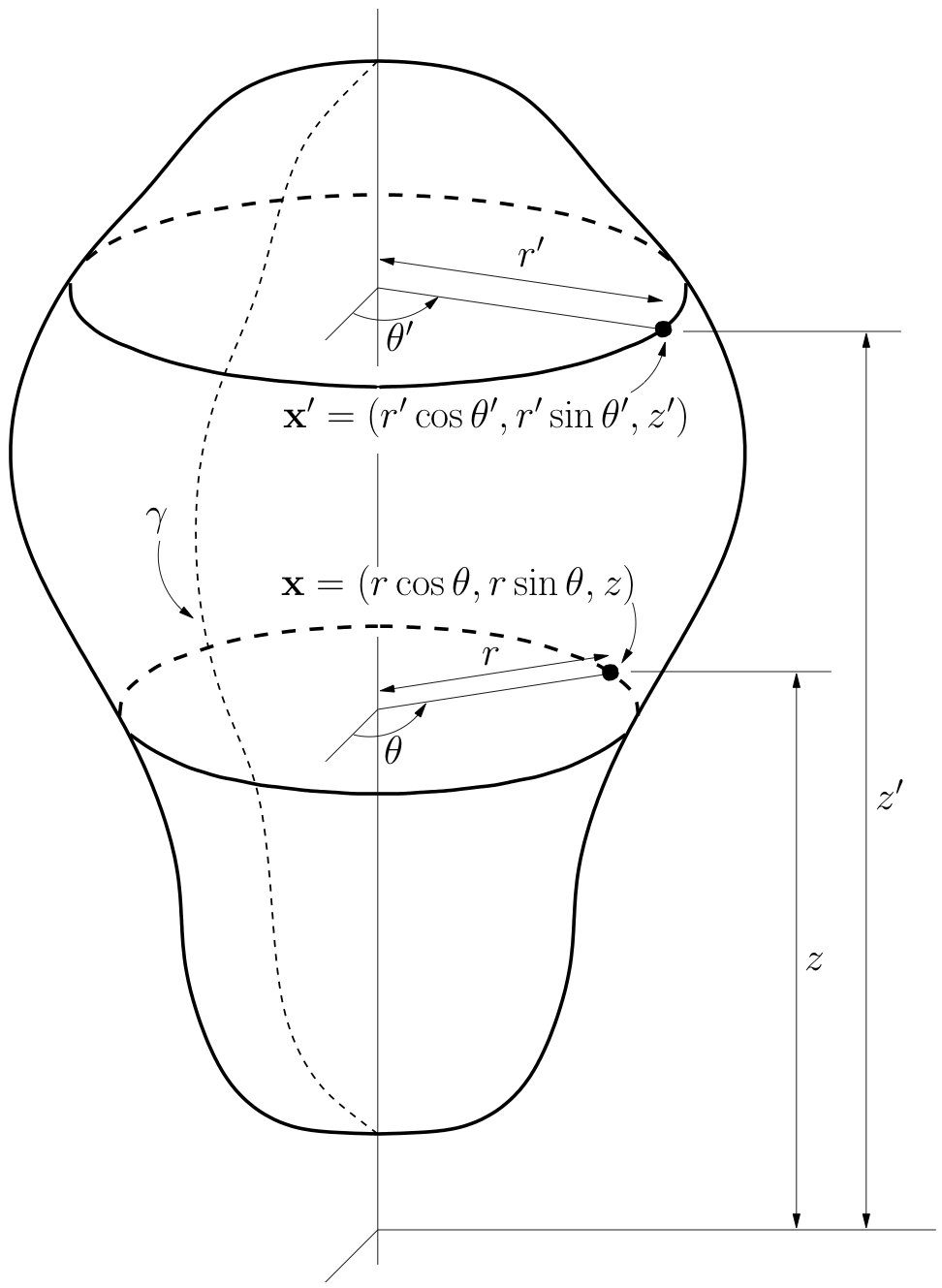}
	\caption{The axisymmetric  domain $\Gamma$ generated by the curve $\gamma$.}
	\label{fig:domain}
\end{figure}

We first consider the case where the scattering surface $\Gamma$ is a single
rotationally symmetric surface. We let $\gamma$ denote a generating curve of $\Gamma$,
and can then view $\Gamma$ as a tensor product between $\gamma$ and the circle $\mathbb{T}$,
so that $\Gamma = \gamma \times \mathbb{T}$, see Figure \ref{fig:domain}. The idea is now to use a composite Gaussian
rule to discretize $\gamma$, and a trapezoidal rule with equispaced nodes to discretize
$\mathbb{T}$, and then take the tensor product between these rules to obtain the  global
rule $\{\pxx_{i},w_{i}\}_{i=1}^{n}$ for $\Gamma$.

\begin{remark}[Convergence order]
Suppose that $\varphi$ is a smooth ($C^{\infty})$ function on $\Gamma$.
Then since $\varphi$ is periodic in the azimuthal direction, the Trapezoidal
rule converges super-algebraically fast. If we use $p$-point Gaussian quadrature
on $r$ intervals to discretize the generating curve $\gamma$, then the error
in (\ref{eq:basicquad}) scales as $(1/r)^{2p-1}$ as $r,p\rightarrow \infty$.
\end{remark}

The technique for constructing a matrix $\mtx{A}$ such that (\ref{eq:nystromquad})
holds is based on the observation that when $\Gamma$ is a rotationally symmetric
surface, the equation (\ref{eq:intHelm2}) is diagonalized by the Fourier transform.
The process is somewhat involved and we will here give only a brief overview of the
key techniques, for details we refer to \cite{2012_martinsson_axisymmetric}.
The first step is to introduce cylindrical coordinates $\pxx = (r,\theta,z)$
with the $z$-axis being the symmetry axis of $\Gamma$,
and let $v_{p}$, $\sigma_{p}$, and $G_{\kappa,p}$ denote the Fourier coefficients of the
functions $v$, $\sigma$, and $G_{\kappa}$:
\begin{align}
\label{eq:forSeries1}
v(\pxx)      &= \sum_{p \in \mathbb{Z}}\frac{e^{ip\theta}}{\sqrt{2\pi}}\,v_{p}(r,z),\\
\label{eq:forSeries2}
\sigma(\pxx) &= \sum_{p \in \mathbb{Z}}\frac{e^{ip\theta}}{\sqrt{2\pi}}\,\sigma_{p}(r,z),\\
\label{eq:forSeries3}
G_{\kappa}(\pxx,\pxx') = G_{\kappa}(\theta-\theta',r,z,r',z') &=
\sum_{p \in \mathbb{Z}}\frac{e^{ip(\theta-\theta')}}{\sqrt{2\pi}}\,G_{\kappa,p}(r,z,r',z').
\end{align}
Then (\ref{eq:intHelm2}) is equivalent to the sequence of equations
\begin{equation}
\label{eq:seqBIE}
\frac{1}{2}\sigma_{p}(\pyy) +
\sqrt{2\pi}\int_{\gamma}G_{\kappa,p}(\pyy,\pyy')\,\sigma_{p}(\pyy')\,dA(\pyy') = -v_{p}(\pyy), \hspace{1em}
\pyy \in \gamma,\ p \in \mathbb{Z}.
\end{equation}
Converting the BIE (\ref{eq:intHelm2}) defined on a surface $\Gamma$ to the sequence
of BIEs (\ref{eq:seqBIE}) defined on the curve $\gamma$ has a crucial advantage in that constructing
high-order Nystr\"om discretizations of BIEs with weakly singular
kernels is well-understood and computationally cheap for curves, but remains a challenge
for surfaces. We use the modified quadrature of \cite{2001_rokhlin_kolm}, as described
in \cite{2012_martinsson_axisymmetric,2013_martinsson_nystrom_final}.

Beyond ease of discretization, the other key benefit of the formulation (\ref{eq:seqBIE}) is that
for each Fourier mode $p$, the coefficient matrix arising from discretization of
(\ref{eq:seqBIE}) is small enough that it can often easily be inverted by brute force.
For instance, for the geometries shown in Figure \ref{fig:domains}, it is sufficient to use at most
a couple of hundred nodes along $\gamma$ to achieve ten digits accuracy. To put it another
way, the Fourier conversion allows to write the matrix $\mtx{A}$ as a product
\begin{equation}
\label{eq:FAF}
\mtx{A} = \mtx{F}^{*}\,\tilde{\mtx{A}}\,\mtx{F}
\end{equation}
where $\mtx{F}$ is the discrete Fourier transform (in the azimuthal variable), and $\tilde{\mtx{A}}$ is a block-diagonal
matrix, where each diagonal block corresponds to one Fourier mode, and is relatively small.
We can pre-compute and store the block diagonal matrix $\tilde{\mtx{A}}^{-1}$, and
then very rapidly apply the inverse
\begin{equation}
\label{eq:singleinverse}
\mtx{A}^{-1} = \mtx{F}^{*}\,\tilde{\mtx{A}}^{-1}\,\mtx{F},
\end{equation}
by using the FFT to apply $\mtx{F}$ and $\mtx{F}^{*}$.

One complication to the procedure outlined in this section is that while the kernel $G_{\kappa}$
in (\ref{eq:intHelm2}) is given by the simple formula (\ref{eq:Gkappa}), the kernels $G_{\kappa,p}$
must be evaluated computationally. Techniques for doing so rapidly have been developed, and
are described in \cite{2012_martinsson_axisymmetric}.

\begin{remark}[Cost of precomputation]
To state the asymptotic cost of the algorithm,
let $N_{\rm G}$ (``G'' for Gaussian) denote the number of points on the generating curve
$\gamma$ of each scatter and let $N_{\rm F}$ (``F'' for Fourier) denote the number of points
used to  discretize $\mathbb{T}$. The total number of degrees of freedom of each scatter is
$n = N_{\rm G}N_{\rm F}$. Under the simplifying assumption that
$N_{\rm G} \sim N_{\rm F}$, the cost of forming the block diagonal
matrix $\tilde{\mtx{A}}$ is $O(n^{3/2}\log{n})$, while the cost of inverting
$\tilde{\mtx{A}}$ is $O(n^2)$, see \cite{2012_martinsson_axisymmetric}.
Applying $\mtx{F}$ and $\mtx{F}^{*}$ is done via the FFT in negligible time.
\end{remark}






\subsection{Multibody scattering}
\label{sec:nystrommulti}
Having described how to discretize the single-body scattering problem in Section \ref{sec:nystromsingle},
we now proceed to the general case of $m$ disjoint scattering surfaces $\Gamma = \cup_{p=1}^{m}\Gamma_{p}$.
We assume that each scatterer is discretized using the tensor product procedure described in Section \ref{sec:nystromsingle}.
For notational simplicity, we assume that each scatterer is discretized using the same $n$ number of nodes,
for a total of $N = mn$ discretization nodes $\{\pxx_{i}\}_{i=1}^{N}$ with associated weights
$\{w_{i}\}_{i=1}^{N}$. We then
seek to construct matrix blocks $\{\mtx{A}_{p,q}\}_{p,q=1}^{m}$ such that the Nystr\"om discretization
of (\ref{eq:intHelm2}) associated with this quadrature rule takes the form
\begin{equation}
\label{eq:blocksystem}
\left[\begin{array}{cccc}
\mtx{A}_{1,1} & \mtx{A}_{1,2} & \cdots & \mtx{A}_{1,m} \\
\mtx{A}_{2,1} & \mtx{A}_{2,2} & \cdots & \mtx{A}_{2,m} \\
\vdots        & \vdots        &        & \vdots        \\
\mtx{A}_{m,1} & \mtx{A}_{m,2} & \cdots & \mtx{A}_{m,m}
\end{array}\right]\,
\left[\begin{array}{c}
\vct{\sigma}_{1} \\ \vct{\sigma}_{2} \\ \vdots \\ \vct{\sigma}_{m}
\end{array}\right]
= -
\left[\begin{array}{c}
\vct{v}_{1} \\ \vct{v}_{2} \\ \vdots \\ \vct{v}_{m}
\end{array}\right],
\end{equation}
where each block $\mtx{A}_{p,q}$ is of size $n\times n$.
The diagonal blocks $\mtx{A}_{p,p}$ are constructed using the technique described in Section \ref{sec:nystromsingle}.
Next observe that in the off-diagonal blocks, the ``naive'' formula (\ref{eq:nystromnaive}) works well since
the kernel $G_{\kappa}(\pxx,\pxx')$ is smooth when $\pxx$ and $\pxx'$ belong to different scatterers.

\begin{remark}
In this paper, we avoid considering the complications of scatterers that touch or are very close.
The procedure described works well as long as the minimal distance between scatterers is not small
compared to the resolution of the quadrature rules used. This means that if two scatterers are moderately
close, high accuracy can be maintained by discretizing these two scatterers more finely.
\end{remark}

\section{A block-diagonal pre-conditioner for the multibody scattering problem}
\label{sec:iteration}

We solve the linear system (\ref{eq:blocksystem}) using the iterative solver GMRES \cite{1986_saad_gmres},
accelerated by a block-diagonal pre-conditioner. To formalize, let us decompose the system matrix as
$$
\mtx{A} = \mtx{D} + \mtx{B},
$$
where
$$
\mtx{D} =
\left[\begin{array}{cccc}
\mtx{A}_{1,1} & \mtx{0} & \mtx{0} & \cdots \\
\mtx{0} & \mtx{A}_{2,2} & \mtx{0} & \cdots \\
\mtx{0} & \mtx{0} & \mtx{A}_{3,3} & \cdots \\
\vdots & \vdots & \vdots &
\end{array}\right]
\qquad\mbox{and}\qquad
\mtx{B} =
\left[\begin{array}{cccc}
\mtx{0} & \mtx{A}_{1,2} & \mtx{A}_{1,3} & \cdots \\
\mtx{A}_{2,1} & \mtx{0} & \mtx{A}_{2,3} & \cdots \\
\mtx{A}_{3,1} & \mtx{A}_{3,2} & \mtx{0} & \cdots \\
\vdots & \vdots & \vdots &
\end{array}\right].
$$
Then we use GMRES to solve the linear system
\begin{equation}
\label{eq:precond}
\vct{\sigma} + \mtx{D}^{-1}\mtx{B}\vct{\sigma} = -\mtx{D}^{-1}\vct{v}.
\end{equation}
We apply the matrix $\mtx{B}$ using the Fast Multipole Method \cite{rokhlin1987,2006_rokhlin_wideband};
specifically the implementation \cite{2011_gimbutas_greengard_FMMLIB3D} by Zydrunas Gimbutas and
Leslie Greengard.

\begin{remark}
\label{remark:costmatvec}
The cost of evaluating the term $\mtx{D}^{-1}\mtx{B}\vct{\sigma}$ in (\ref{eq:precond}) consists of two parts: applying $\mtx{B}$ to vector $\vct{\sigma}$ via FMM costs $O(mn)$ operations and applying the block-diagonal pre-conditioner costs $O(mn^{3/2})$
operations. Observe that the matrix $\mtx{D}^{-1}$ can be precomputed since each matrix $\mtx{A}_{p,p}^{-1}$ is itself block-diagonal in the local Fourier basis, cf.~formula (\ref{eq:singleinverse}). Applying $\mtx{A}_{p,p}^{-1}$ to a vector $\vct{w} \in \mathbb{C}^{n}$ is executed as follows:
(1) form $\mtx{F}\vct{w}$ using the FFT at cost $O(n\log{n})$,
(2) for each Fourier mode apply $\mtx{D}^{-1}$ to $\mtx{F}\vct{w}$ at cost $O(n^{3/2})$,
and (3) use the FFT to apply $\mtx{F}^{*}$ to $\mtx{D}^{-1}\mtx{F}\vct{w}$.
\end{remark}

\section{Accelerated multibody scattering}
\label{sec:accelerate}

In situations where the scatterers are not tightly packed, it is often possible to
substantially reduce the size of the linear system (\ref{eq:precond}) before applying
an iterative solver. We use a technique that was introduced in \cite{lowrank} for
problems in two dimensions, which exploits that when the scatterers are somewhat
separated, the off-diagonal blocks $\mtx{A}_{p,q}$ are typically rank deficient.
Specifically, we assume that for some finite precision $\varepsilon$ (say $\varepsilon =10^{-10}$),
each such block admits a factorization
\begin{equation}
\label{eq:ip0}
\begin{array}{ccccccc}
\mtx{A}_{p,q} &=& \mtx{U}_{p} & \tilde{\mtx{A}}_{p,q} & \mtx{V}_{q}^{*} & + & \mtx{R}_{p,q}\\
n\times n && n\times k & k\times k & k\times n && n\times n
\end{array}
\end{equation}
where $n$ is the number of nodes originally used to discretize a single scatterer, and $k$
is the numerical rank of the factorization. The remainder term $\mtx{R}_{p,q}$ satisfies
$||\mtx{R}_{p,q}|| \leq \varepsilon$ in some suitable matrix norm (we typically use the Frobenius norm since
it is simple to compute).

Now write the linear system (\ref{eq:precond}) in block form as
\begin{equation}
\label{eq:ip1}
\vct{\sigma}_{p}
+
\sum_{q \neq p}
\mtx{A}_{p,p}^{-1}\mtx{A}_{p,q}\vct{\sigma}_{q}
=
-\mtx{A}_{p,p}^{-1}\vct{v}_{p},\qquad p = 1,\,2,\,3,\,\dots,\,m.
\end{equation}
We left multiply (\ref{eq:ip1}) by $\mtx{V}_{p}^{*}$, and insert the factorization (\ref{eq:ip0}) to obtain
\begin{equation}
\label{eq:ip2}
\mtx{V}_{p}^{*}\vct{\sigma}_{p}
+
\sum_{q \neq p}
\mtx{V}_{p}^{*}\mtx{A}_{p,p}^{-1}\mtx{U}_{p}\tilde{\mtx{A}}_{p,q}\mtx{V}_{q}^{*}\vct{\sigma}_{q}
=
-\mtx{V}_{p}^{*}\mtx{A}_{p,p}^{-1}\vct{v}_{p},\qquad p = 1,\,2,\,3,\,\dots,\,m.
\end{equation}
We now define quantities
$\{\tilde{\vct{\sigma}}_{p}\}_{p=1}^{m}$,
$\{\tilde{\vct{v}}_{p}\}_{p=1}^{m}$,
and
$\{\tilde{\mtx{S}}_{p}\}_{p=1}^{m}$
via
\begin{equation}
\label{eq:ip3}
\tilde{\vct{\sigma}}_{p} = \mtx{V}_{p}^{*}\vct{\sigma}_{p},
\qquad
\tilde{\vct{v}}_{p} = \mtx{V}_{p}^{*}\mtx{A}_{p,p}^{-1}\vct{v}_{p}
\qquad
\mtx{S}_{p,p} = \mtx{V}_{p}^{*}\mtx{A}_{p,p}^{-1}\mtx{U}_{p},\qquad\mbox{for}\ p = 1,\,2,\,3,\,\dots,\,m.
\end{equation}
Then the system (\ref{eq:ip2}) can be written
\begin{equation}
\label{eq:ip4}
\tilde{\vct{\sigma}}_{p}
+
\sum_{q \neq p}
\mtx{S}_{p}\tilde{\mtx{A}}_{p,q}\tilde{\vct{\sigma}}_{q}
=
-\tilde{\vct{v}}_{p},\qquad p = 1,\,2,\,3,\,\dots,\,m.
\end{equation}
To write (\ref{eq:ip4}) in block form, we introduce matrices
\begin{equation}
\label{eq:ip5}
\mtx{S} =
\left[\begin{array}{cccc}
\mtx{S}_{1} & \mtx{0} & \mtx{0} & \cdots \\
\mtx{0} & \mtx{S}_{2} & \mtx{0} & \cdots \\
\mtx{0} & \mtx{0} & \mtx{S}_{3} & \cdots \\
\vdots & \vdots & \vdots &
\end{array}\right]
\qquad\mbox{and}\qquad
\tilde{\mtx{B}} =
\left[\begin{array}{cccc}
\mtx{0} & \tilde{\mtx{A}}_{1,2} & \tilde{\mtx{A}}_{1,3} & \cdots \\
\tilde{\mtx{A}}_{2,1} & \mtx{0} & \tilde{\mtx{A}}_{2,3} & \cdots \\
\tilde{\mtx{A}}_{3,1} & \tilde{\mtx{A}}_{3,2} & \mtx{0} & \cdots \\
\vdots & \vdots & \vdots &
\end{array}\right],
\end{equation}
whence equation (\ref{eq:ip4}) takes the form, cf.~(\ref{eq:precond}),
\begin{equation}
\label{eq:ip6}
\tilde{\vct{\sigma}} + \mtx{S}\tilde{\mtx{B}}\tilde{\vct{\sigma}} = -\tilde{\vct{v}}.
\end{equation}

The process of first forming the linear system (\ref{eq:ip6}), and then solving it
using GMRES is very computationally efficient when the following techniques are used:
\begin{itemize}
\item The matrices $\{\mtx{U}_{p},\mtx{V}_{p}\}_{p=1}^{m}$ in the factorizations (\ref{eq:ip0})
can be computed via a purely local procedure in $O(n^{2}k)$ operations, independent of the
number of scatterers $m$. The idea is to use representation techniques from scattering theory
to construct a local basis for all possible incoming harmonic fields (to within
precision $\varepsilon$), see \cite[Sec.~5.1]{2009_martinsson_ACTA} or \cite[Sec.~6.2]{2012_martinsson_FDS_survey}.
\item In constructing the factorization (\ref{eq:ip0}), the so called \textit{interpolatory decomposition}
\cite{lowrank} should be used. Then each matrix $\mtx{U}_{p}$ and each matrix $\mtx{V}_{p}$ contains
the $k\times k$ identity matrix $\mtx{I}_{k}$. Specifically, there exists for each $k$ an index
vector $\tilde{I}_{p} \subset \{1,2,\dots,n\}$ such that $\mtx{U}(\tilde{I}_{p},:) = \mtx{V}(\tilde{I}_{p},:) = \mtx{I}_{k}$.
Then each off-diagonal block $\tilde{\mtx{A}}_{p,q}$ is given as a \textit{submatrix}
$\tilde{\mtx{A}}_{p,q} = \mtx{A}_{p,q}(\tilde{I}_{p},\tilde{I}_{q})$. In consequence,
the matrix $\tilde{\mtx{B}}$ is a sub-matrix of $\mtx{B}$ and can be rapidly applied using the FMM
in $O(mk)$ operations.
\item In evaluating the formula $\mtx{S}_{p,p} = \mtx{V}_{p}^{*}\mtx{A}_{p,p}^{-1}\mtx{U}_{p}$,
we exploit that $\mtx{A}_{p,p}^{-1}$ can be applied rapidly in Fourier space, cf.~(\ref{eq:singleinverse}),
to reduce the complexity of this step from $O(n^{3})$ to $O(n^{3/2}k)$ if $\mtx{A}_{p,p}^{-1}$ was precomputed and stored and to $O(n^{2}k)$ if $\mtx{A}_{p,p}^{-1}$ is computed at this step.
\end{itemize}

\begin{remark}
Efficient techniques for computing interpolative decompositions are described in
\cite{lowrank}. More recently, techniques based on randomized sampling have proven
to be highly efficient on modern computing platforms, in particular for problems
in potential theory where the low-rank matrices to be approximated have very rapidly
decaying singular values. We use the specific technique described in \cite{MR2806637}.
\end{remark}

\section{Numerical examples}
\label{sec:num}

This section describes numerical experiments to assess the performance of the numerical
scheme outlined in previous sections. All the experiments are carried out on a personal
work-station with an Intel Xeon E-1660 3.3GHz 6-core CPU, and 128GB of RAM. The
experiments explore (1) the accuracy of the algorithm, (2) the computational cost,
(3) the performance of block-diagonal pre-conditioner and (4) the performance of
the acceleration scheme when scatterers are separated suitably.
In all the experiments below, we measure accuracy against a known analytic solution $u_{\rm exact}$.
This solution is generated by randomly placing one point source inside each scatterer, and
then solving (\ref{eq:basic}) with the Dirichlet data $v$ set to equal the field generated by these
radiating sources.
Let $\vct{u}_{\textrm{exact}}$ and $\vct{u}_{\textrm{approx}}$ denote the vectors holding the
exact and the computed solutions at a set of 10 randomly chosen target points, located a few
wave-lengths outside of the domain. The relative error, measured in $\ell^{\infty}$-norm,
is then given by
$$
E^{\textrm{rel}}_{\infty} = \frac{||\vct{u}_{\textrm{approx}} - \vct{u}_{\textrm{exact}}||_{\infty}}{ || \vct{u}_{\textrm{exact}} ||_{\infty}}.
$$
In addition to $E^{\rm rel}_{\infty}$, we report:\\
\begin{tabular}{l l}
$n$ &  number of nodes discretizing each body (in form of $n = N_{\rm G} \times N_{\rm F}$)\\
$N$ & total degree of freedom $N = m\times n$, where $m$ is the number of scatterers\\
$N_{\textrm{compressed}}$ & number of skeleton points after applying the compression scheme\\
$T_{\textrm{pre}}$ & time (in seconds) of precomputation\\
$T_{\textrm{solve}}$ & total time to solve for the surface charges $\sigma$ via GMRES\\
$T_{\textrm{compress}}$ & time to do compression in the accelerated scheme \\
$I$ & number of GMRES iterations required to reduce the residual to $10^{-9}$.
\end{tabular}

\begin{figure}[h!]
\begin{center}
\begin{minipage}{0.31\linewidth}
\begin{flushleft}
\hspace{-9mm}
\includegraphics[height=.93\linewidth]{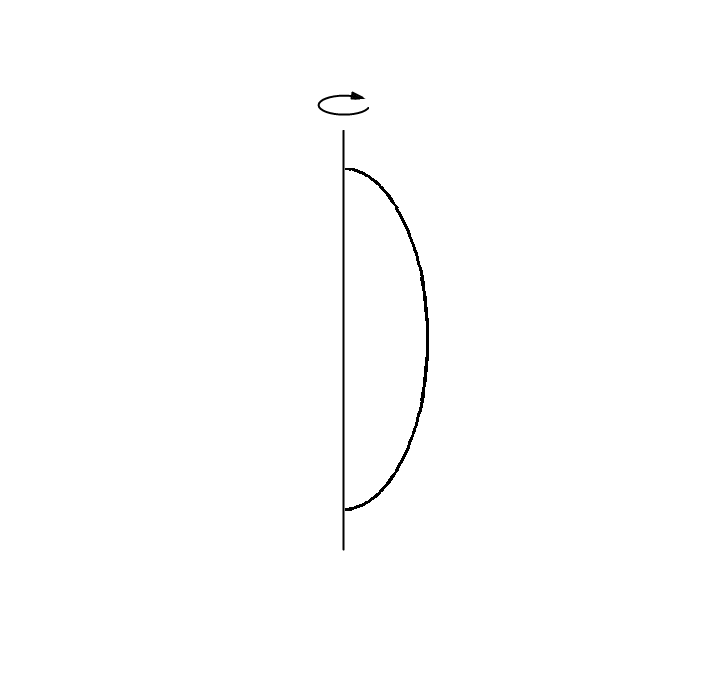}
\end{flushleft}
\end{minipage}
\hspace{2mm}
\begin{minipage}{0.31\linewidth}
\begin{center}
\includegraphics[height =.8\linewidth]{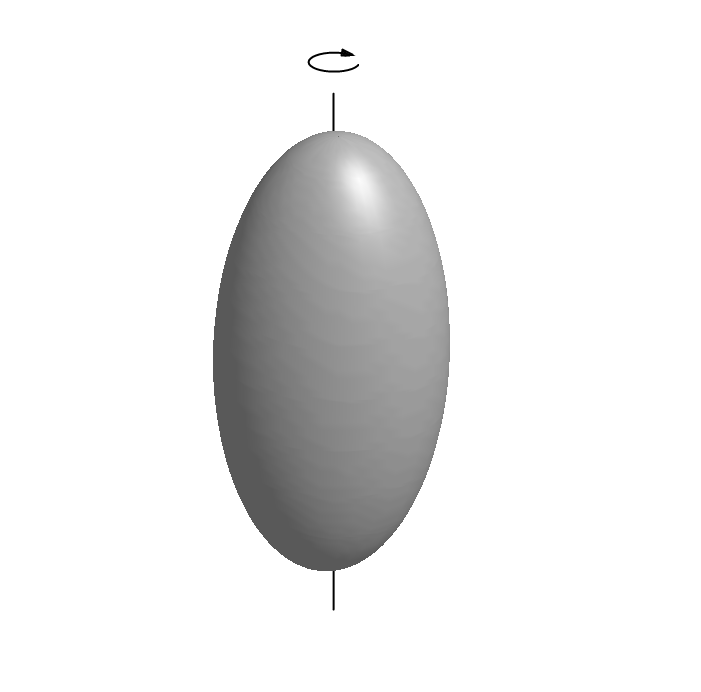}
\end{center}
\end{minipage} \\
\begin{minipage}{1\linewidth}\begin{center} (a) \end{center}\end{minipage} \\
\vspace{4mm}
\begin{minipage}{0.31\linewidth}
\begin{flushleft}
\hspace{10mm}
 \includegraphics[height=.6\linewidth]{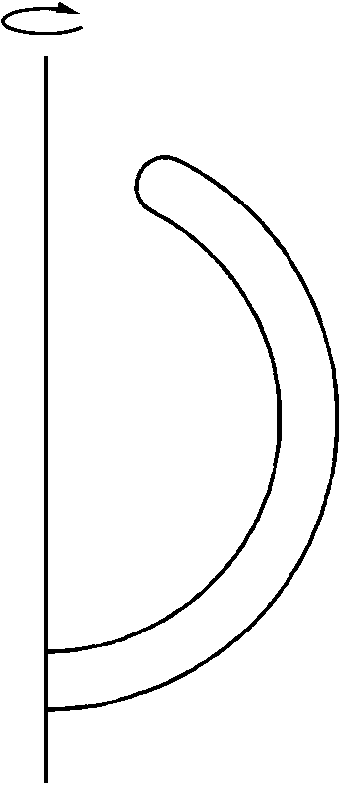}
\end{flushleft}
\end{minipage}
\begin{minipage}{0.31\linewidth}
\begin{center}
\includegraphics[height =.5\linewidth]{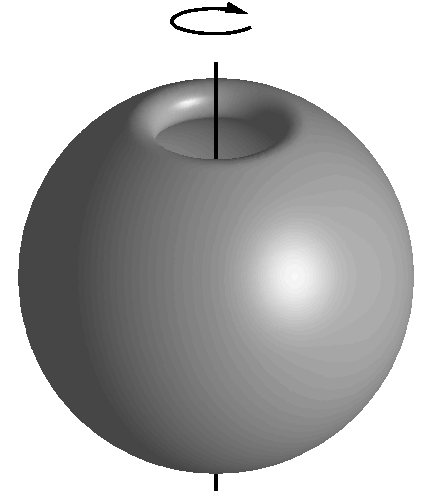}
\end{center}
\end{minipage} \\
\begin{minipage}{1\linewidth}\begin{center} (b) \end{center}\end{minipage} \\
\vspace{4mm}
\begin{minipage}{0.31\linewidth}
\begin{flushleft}
\hspace{0mm}
\includegraphics[height=.75\linewidth]{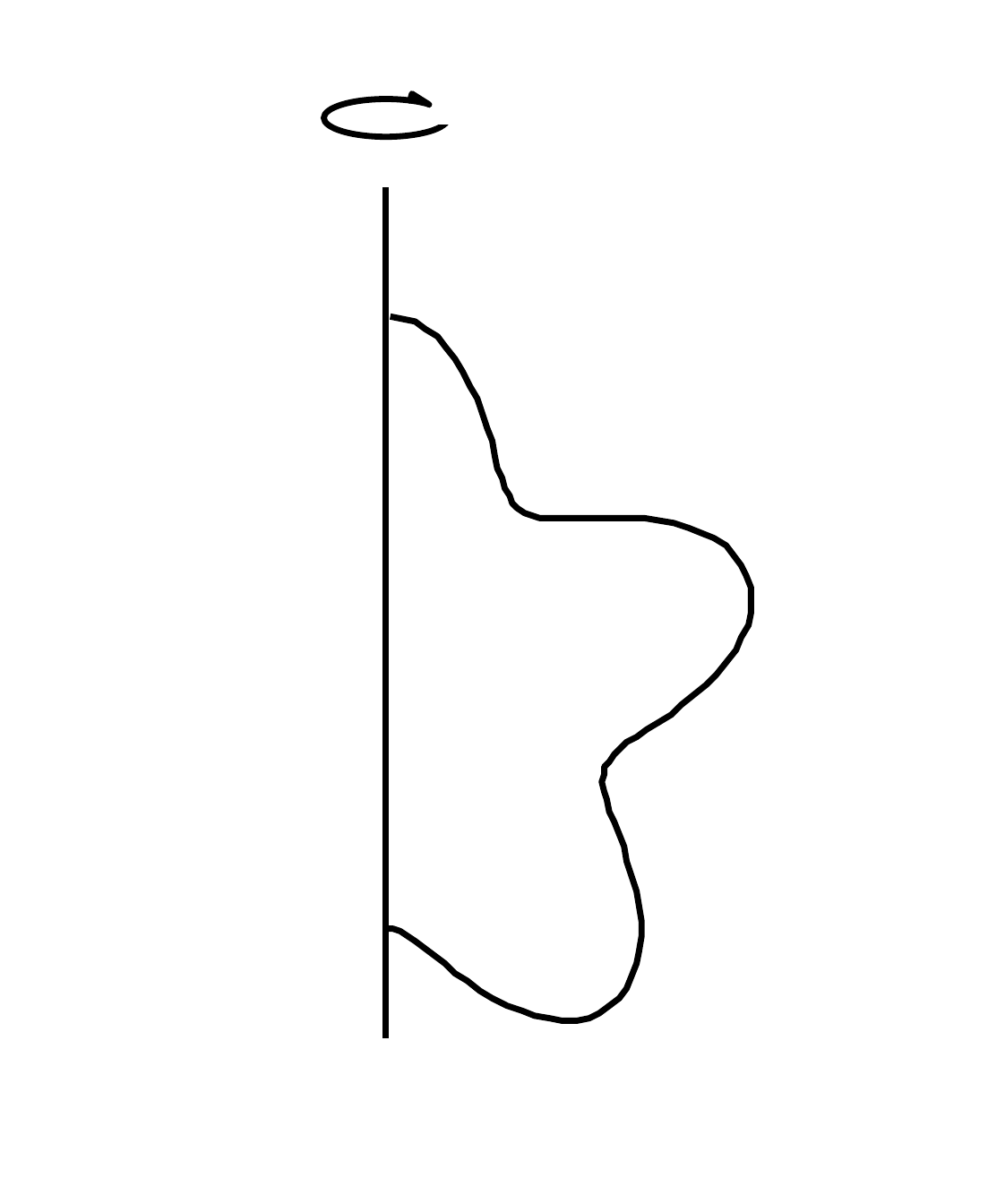}
\end{flushleft}
\end{minipage}
\begin{minipage}{0.31\linewidth}
\begin{center}
\includegraphics[height =.65\linewidth]{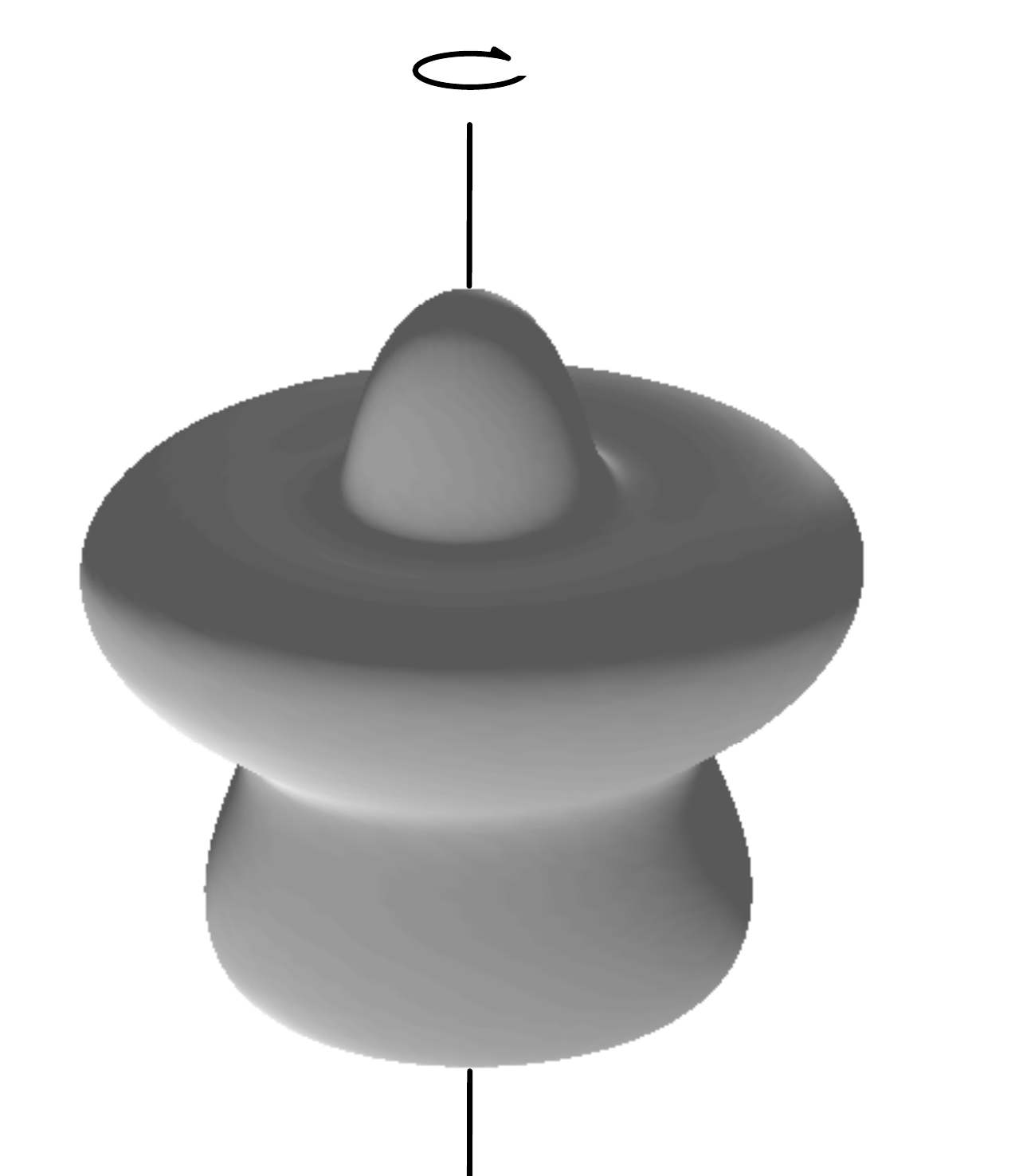}
\end{center}
\end{minipage} \\
\begin{minipage}{1\linewidth}\begin{center} (c) \end{center}\end{minipage} \\
\caption{Domains used in numerical examples. All items are rotated
about their symmetry axis.  (a) An ellipsoid.  (b) A bowl-shaped cavity.  (c) A starfish-shaped cavity.}
\label{fig:domains}
\end{center}
\end{figure}

All the numerical experiments in this section are executed on domains composed of the
three sample scatterers shown in Figure \ref{fig:domains}.

\subsection{Laplace's equation}
\label{sec:numLap}

We first solve the Laplace equation exterior to the domain shown in Figure \ref{fig:ellp_125} and \ref{fig:cavity_8}(a) (Example 1 and 2 respectively). A combination of the single and double layer kernels is chosen to represent the potential outside the domain.  The integral equation to be solved is
$$
	\frac{1}{2}\sigma(\xb)+\int_{\Gamma}\frac{1}{4\pi}\(\frac{1}{|\xb-\xb'|}+\frac{\nb(\xb')\cdot (\xb-\xb')}{|\xb-\xb'|^3}\)\,\sigma(\xb')\,dA(\xb') = f(\xb), \hspace{1em} \xb \in \Gamma.
$$

\subsubsection{Example 1}
\begin{figure}[htbp]
\begin{center}
\includegraphics[height =.35\linewidth]{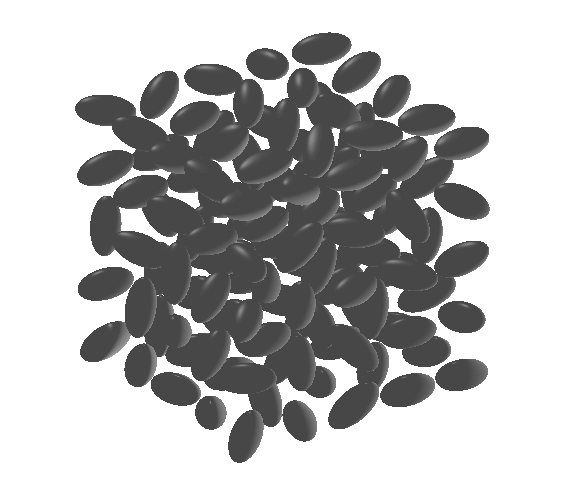}
\end{center}
\caption{Domain contains 125 randomly oriented ellipsoids. Distance between any two ellipsoids is $0.05$.}
\label{fig:ellp_125}
\end{figure}

This example solves the exterior Laplace equation on the domain depicted in Figure \ref{fig:ellp_125}. The domain consists of 125 ellipsoids contained in the box $[0,10.2]^3$, where each ellipse has a major axis of length 2 and a minor axis of length 1. The minimal distance between any two ellipsoids is 0.05. We did not apply the compression technique since the scatterers are packed tightly. We compare the performance of the algorithm with and without using block-diagonal pre-conditioner in Table \ref{tab:lap_125ellp} and find that for this example, the pre-conditioning does not make any real difference. The scheme quickly reaches 9 digits of accuracy with $10\,100$ discretization nodes per scatterer, with an overall solve time of
about 40 minutes.

\begin{table}[htbp]
\centering
\begin{tabular}{|r|r|c|c|c|c|}
\hline
\multirow{2}{*}{$N$} & \multirow{2}{*}{$n$}  &  \multirow{2}{*}{$T_{\textrm{pre}}$} & $I$ & $T_{\textrm{solve}} $  & \multirow{2}{*}{$E^{\textrm{rel}}_{\infty}$}\\  & & & (precond  /no precond ) & (precond /no precond) &  \\
\hline
156\,250 &$50\times   25$ & 1.09e+00  & 31 /33 & 3.16e+02 /3.29e+02 & 9.731e-05\\
312\,500 & $100\times 25$& 3.44e+00  & 31 /33 & 6.84e+02 /6.82e+02 &  9.203e-05\\
625\,000 & $200\times 25$& 1.29e+01  & 31 /34 & 1.10e+03 /1.18e+03 &  9.814e-05\\
\hline
318\,750 & $50\times  51$ & 1.53e+00  & 31 /33 & 6.29e+02 /7.44e+02 & 1.571e-06\\
637\,500 & $100\times 51$& 4.36e+00 & 31 /34 & 1.18e+03 /1.23e+03 &  1.529e-06\\
1\,275\,000&$200\times 51$& 1.36e+01  & 32 /34 & 2.70e+03 /2.53e+04 & 1.711e-06\\
\hline
631\,250 & $50\times  101$& 2.44e+00  & 31 /34 & 1.11e+03 /1.22e+03 & 2.165e-08\\
1\,262\,500&$100\times 101$& 6.11e+00  & 32 /34& 2.45e+03 /2.60e+03&  1.182e-09\\
\hline
\end{tabular}
\vspace{2mm}
\caption{Example 1: exterior Laplace problem solved on the domain in Figure \ref{fig:ellp_125}.}
\label{tab:lap_125ellp}
\end{table}

\subsubsection{Example 2}
\begin{figure}[htbp]
\begin{center}
\begin{tabular}{ccc}
\includegraphics[height=.45\linewidth]{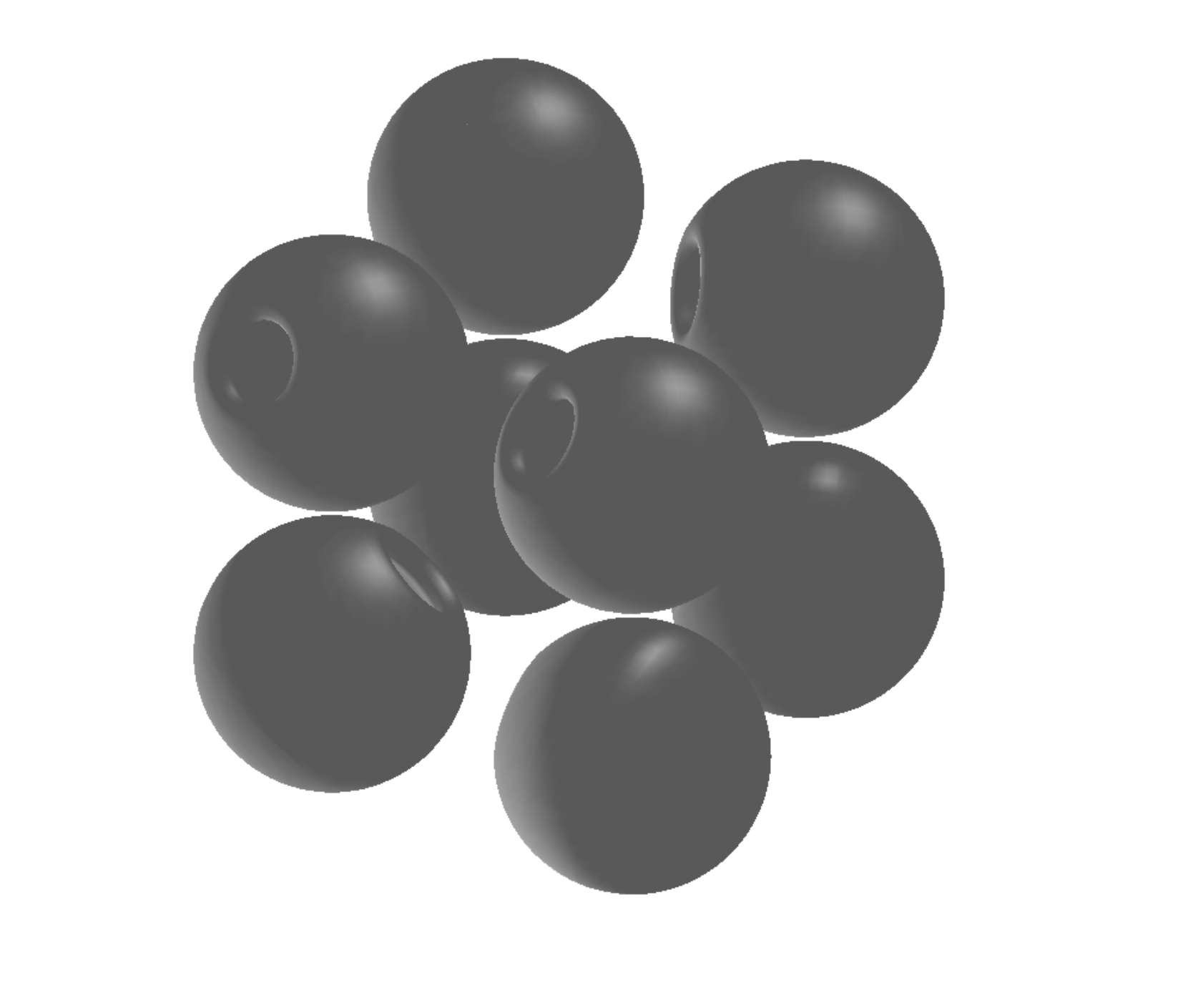}
&\mbox{}\hspace{5mm}\mbox{}&
\includegraphics[width=0.45\textwidth]{figs/multi_ellipses8.pdf}\\
(a) && (b)
\end{tabular}
\end{center}
\caption{(a) Domain contains 8 bowl-shaped cavities. Distance between any two cavities is $0.5$.
(b) Domain contains 8 randomly oriented ellipsoids. Distance between any two ellipsoids is $0.05$.}
\label{fig:cavity_8}
\end{figure}

This time the domain consists of 8 bowl-shaped cavities contained in the box $[0,4.1]^3$ in Figure \ref{fig:cavity_8}(a). The minimal distance between any two cavities is 0.5. Results are shown in Table \ref{tab:lap_cavity}. The scheme achieves 8 digits of accuracy with $400$ discretization nodes on the generating curve and $201$ Fourier modes. Again, the pre-conditioning is superfluous.

\begin{table}[htbp]
\centering
\begin{tabular}{|r|r|c|c|c|c|}
\hline
\multirow{2}{*}{$N$} & \multirow{2}{*}{$n$}  &  \multirow{2}{*}{$T_{\textrm{pre}}$} & $I$ & $T_{\textrm{solve}} $  & \multirow{2}{*}{$E^{\textrm{rel}}_{\infty}$}\\  & & & (precond  /no precond ) & (precond /no precond) &  \\
\hline
20\,400&50$\times$51  & 2.09e-01 &398 /402 & 4.65e+02 /6.05e+02 &1.251e-04\\
40\,800&100$\times$51& 4.55e-01 & 20 /23   & 4.94e+01 /6.09e+01   &3.909e-05\\
81\,600&200$\times$51& 9.83e-01 & 20 /23   & 1.05e+02 /1.14e+02   &3.164e-05\\
\hline
40\,400&50$\times$101    & 2.25e-01 & 20 /23    & 4.72e+01 /6.17e+01  & 5.850e-05\\
80\,800&100$\times$101  & 4.49e-01 & 20 /23    & 9.50e+01 /1.13e+02  & 1.627e-05\\
161\,600&200$\times$101&1.35e+00 & 20 /24     & 2.05e+02 /2.39e+02  & 6.825e-06\\
\hline
80\,400 & $50 \times 201$     &  2.93e-01 & 20 /23 & 9.13e+01 /1.12e+02 & 5.704e-05 \\
160\,800 & $100 \times 201$ & 7.05e-01  & 20 /24 & 1.96e+02 /2.40e+02 & 8.000e-06 \\
321\,600 & $200 \times 201$ & 1.97e+00 & 20 /24 & 4.43e+02 /5.25e+02 &1.931e-07\\
643\,200 & $400 \times 201$ & 5.78e+00 & 21 /24 & 7.68e+02 /8.19e+02 & 1.726e-08 \\
\hline
\end{tabular}
\vspace{2mm}
\caption{Example 2: exterior Laplace problem solved on the domain in Figure \ref{fig:cavity_8}(a).}
\label{tab:lap_cavity}
\end{table}

\begin{remark}
All examples described in this section involve geometries where all
the scatterers are copies of the basic shapes shown in Figure \ref{fig:domains}.
In our experience, this restriction on the geometry does not in any
way change the overall accuracy or efficiency of the solver.
The only advantage we benefit from is that the pre-computation
gets faster, as only a small number of scattering matrices need to
be pre-computed. However, it is clear from the numbers given that
even for a fully general geometry (without repetitions), the
pre-computation time would be dominated by the time required for
the FMM.
\end{remark}

\subsection{Helmholtz Equation}
\label{sec:numHelm}
We now consider the exterior Helmholtz problem (\ref{eq:basic}).
We represent the potential by a combination of the single and double layer kernels,
see (\ref{eq:Gkappa}), and end up with the ``combined field'' integral equation (\ref{eq:intHelm2}).

\subsubsection{Example 3} The domain in this experiment contains 8 ellipsoids in the box $[0,4.05]^3$, whose minimal distance between any two is 0.05. The wavelength is $10\pi$ so that the scatterers are approximately 10 wavelength in size and the whole region is about $20\times20\times20$ wavelengths in size. Results are presented in Table \ref{tab:helm_8ellp_10}. We also compare the results without using block-diagonal pre-conditioner in the same table. Around twice of the iteration numbers are required resulting in twice of the computation time. Table \ref{tab:helm_8ellp_20} reports the results from an analogous experiment, but now the wavenumber increases such that each scatterer contains 20 wavelengths.

\begin{table}[htbp]
\centering
\begin{tabular}{|r|r|c|c|c|c|}
\hline
\multirow{2}{*}{$N$} & \multirow{2}{*}{$n$}  &  \multirow{2}{*}{$T_{\textrm{pre}}$} & $I$ & $T_{\textrm{solve}} $  & \multirow{2}{*}{$E^{\textrm{rel}}_{\infty}$}\\  & & & (precond  /no precond ) & (precond /no precond) &  \\
\hline
20\,400&$50\times   51$& 1.58e-01     & 35 /67   & 7.71e+02 /1.56e+03  & 1.364e-03\\
40\,800&$100\times 51$& 4.20e-01     & 36 /67  & 1.75e+03 /3.43e+03  &1.183e-03\\
81\,600&$200\times 51$& 1.26e-01     & 36 /68  & 3.52e+03 /6.85e+03   &1.639e-04\\
\hline
40\,400& $50\times   101$ & 2.64e-01 & 36 /68 & 1.71e+03 /3.35e+03  &1.312e-03\\
80\,800& $100\times 101$ & 6.05e-01 & 36 /68 & 3.45e+03 /6.76e+03  &1.839e-06\\
161\,600&$200\times101$ & 1.87e+00 & 37 /69 & 6.18e+03 /1.19e+04  &5.126e-08\\
\hline
80\,400&  $50\times   201$ & 4.61e-01 & 36 /69 & 3.40e+03 /6.70e+03  &1.312e-03\\
160\,800&$100\times 201$ & 1.09e+00 & 37 /69 & 6.07e+03 /1.18e+04 &1.851e-06\\
321\,600&$200\times 201$ & 3.11e+00 & 37 /69 & 1.20e+04 /1.97e+04  & 1.039e-09 \\
\hline
\end{tabular}
\vspace{2mm}
\caption{Example 3: exterior Helmholtz problem solved on the domain in Figure \ref{fig:cavity_8}(b). Each ellipsoid is 10 wavelength in diameter.}
\label{tab:helm_8ellp_10}
\end{table}

\begin{table}[htbp]
\centering
\begin{tabular}{|r|r|c|c|c|c|}
\hline
\multirow{2}{*}{$N$} & \multirow{2}{*}{$n$}  &  \multirow{2}{*}{$T_{\textrm{pre}}$} & $I$ & $T_{\textrm{solve}} $  & \multirow{2}{*}{$E^{\textrm{rel}}_{\infty}$}\\  & & & (precond  /no precond ) & (precond /no precond) &  \\
\hline
20\,400&$50\times   51$&2.03e-01       &  58 /119   & 3.59e+03 /8.10e+03  & 4.362e+00\\
40\,800&$100\times 51$&4.44e-01        & 39 /102   & 3.98e+03 /1.11e+04 & 1.071e+00\\
81\,600&$200\times 51$&1.36e+00      &  39 /106  & 6.72e+03 / 1.92e+04 & 1.008e+00 \\
\hline
40\,400& $50\times   101$&2.78e-01 &  54 /94   & 5.43e+03 /1.02e+04   &5.039e+00\\
80\,800& $100\times 101$&6.18e-01 &  36 /82  & 6.11e+03  /1.46e+04   &8.919e-04\\
161\,600&$200\times101$&1.93e+00 & 36 /83  & 9.44e+03 /2.32e+04   &5.129e-07\\
\hline
80\,400&  $50\times   201$ & 4.28e-01   &  55 /95 &  9.19e+03 /6.41e+04   & 5.031e+00 \\
160\,800&$100\times 201$ &1.07e+00   &  36 /83 &  9.49e+03 /2.31e+04   &8.916e-04\\
321\,600&$200\times 201$ &3.10e+00   &  37 /83 &  1.45e+04 /3.57e+04 &8.781e-09   \\
\hline
\end{tabular}
\vspace{2mm}
\caption{Example 3: exterior Helmholtz problem solved on the domain in Figure \ref{fig:cavity_8}(b). Each ellipsoid is 20 wavelength in diameter.}
\label{tab:helm_8ellp_20}
\end{table}

\subsubsection{Example 4}
This example solves the exterior Helmholtz problem on the cavity domain in Figure \ref{fig:cavity_8}(a).
Tables \ref{tab:helm_cavity_2} and \ref{tab:helm_cavity_5} show the results from experiments involving
cavities of diameters 2 and 5 wavelengths, respectively. In this case, computing the actual scattering
matrix for each scatterer was \textit{essential}, without using these to pre-condition the problem, we
did not observe any convergence in GMRES.

\begin{table}[htbp]
\centering
\begin{tabular}{|r|r|c|c|c|c|}
\hline
\multirow{2}{*}{$N$} & \multirow{2}{*}{$n$}  &  \multirow{2}{*}{$T_{\textrm{pre}}$} & $I$ & $T_{\textrm{solve}} $  & \multirow{2}{*}{$E^{\textrm{rel}}_{\infty}$}\\  & & & (precond  /no precond ) & (precond /no precond) &  \\
\hline
40\,800&$100\times 51$& 4.29e-01    & 59 /181  & 2.17e+03 /6.73e+03  &1.127e-02\\
81\,600&$200\times 51$& 1.28e+00    & 60 / --  & 4.23e+03 / --  &1.131e-02\\
\hline
80\,800& $100\times 101$ & 6.83e-01  & 60 / -- & 4.18e+03 / --  &3.953e-03\\
161\,600&$200\times101$ & 1.90e+00  & 60 / -- & 8.93e+03 / --  &3.802e-04\\
323\,200& $400\times101$ & 6.07e+00 & 61 / --& 1.91e+04 / --  & 3.813e-04\\
\hline
160\,800&$100\times 201$ & 1.09e+00 & 60 / --  & 8.35e+03 / --  & 4.788e-05\\
321\,600&$200\times 201$ & 3.07e+00 & 61 / -- & 1.88e+04 / --  & 5.488e-06 \\
 643\,200&$400\times201$ & 9.61e+00 & 61 / -- & 4.03e+04 / --  &8.713e-08 \\
\hline
\end{tabular}
\vspace{2mm}
\caption{Example 4: exterior Helmholtz problem solved on the domain in Figure \ref{fig:cavity_8}(a). Each cavity is 2 wavelength in diameter.}
\label{tab:helm_cavity_2}
\end{table}

\begin{table}[htbp]
\centering
\begin{tabular}{|r|r|c|c|c|c|}
\hline
\multirow{2}{*}{$N$} & \multirow{2}{*}{$n$}  &  \multirow{2}{*}{$T_{\textrm{pre}}$} & $I$ & $T_{\textrm{solve}} $  & \multirow{2}{*}{$E^{\textrm{rel}}_{\infty}$}\\  & & & (precond  /no precond ) & (precond /no precond) &  \\
\hline

80\,800& $100\times 101$ & 6.54e-01  & 62 /304 & 5.17e+03 / 2.64e+04 &1.555e-03\\
161\,600&$200\times101$ & 1.82e+00  & 63 / -- & 9.88e+03 / -- &1.518e-04\\
323\,200& $400\times101$ & 6.46e+00 & 64 / -- & 2.19e+04 / -- & 3.813e-04\\
\hline
160\,800&$100\times 201$ & 1.09e+00 & 63 / -- & 9.95e+03 / -- & 1.861e-03\\
321\,600&$200\times 201$ & 3.00e+00 & 64 / -- & 2.19e+04 /  -- & 2.235e-05 \\
 643\,200&$400\times201$ & 1.09e+01 & 64 / -- & 4.11e+04 / --  &8.145e-06 \\
\hline
641\,600&$200\times 401$ & 5.02e+00 & 64 / -- & 4.07e+04 / -- &2.485e-05 \\
1\,283\,200&$400\times401$ & 1.98e+01  & 65 / -- & 9.75e+04 / -- &6.884e-07 \\
\hline
\end{tabular}
\vspace{2mm}
\caption{Example 4: exterior Helmholtz problem solved on the domain in Figure \ref{fig:cavity_8}(a). Each cavity is 5 wavelength in diameter.}
\label{tab:helm_cavity_5}
\end{table}


\subsection{Accelerated scheme}
\label{sec:numCompress}

\subsubsection{Example 5}
\begin{figure}[htbp]
\centering
\includegraphics[width=0.8\textwidth]{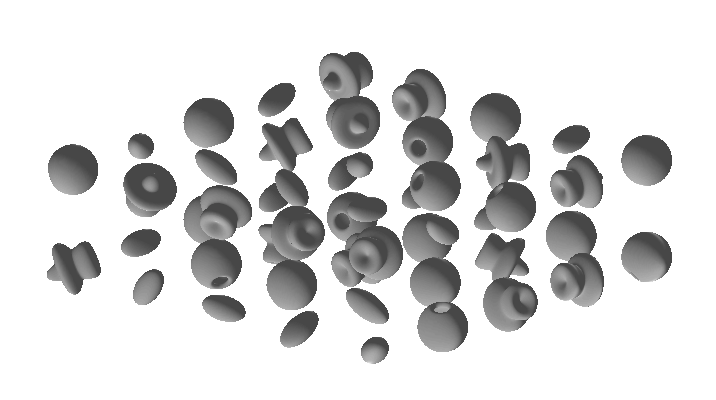}
\caption{Domain contains $50$ randomly oriented scatters.}
\label{fig:diff_scatter}
\end{figure}

We apply the accelerated scheme in Section \ref{sec:accelerate} to solve the Laplace's equation on the
domain exterior to the bodies depicted in Figure \ref{fig:diff_scatter}. This geometry contains 50 different
shaped scatterers (ellipsoids, bowls, and rotated ``starfish'') and is contained in the box
$[0,18]\times[0,18]\times[0,6]$. The minimal
distance between any two bodies is $4.0$. We compare the results shown in Table \ref{tab:lap_comp} with the
ones without applying accelerated scheme shown in Table \ref{tab:lap_nocomp}.
We see that the compression did not substantially alter either the convergence speed of GMRES, or the final
accuracy, and therefore led to improvements in the solve time of between one and two orders of magnitude.

\begin{figure}
\begin{tabular}{ccccc}
\includegraphics[width=30mm]{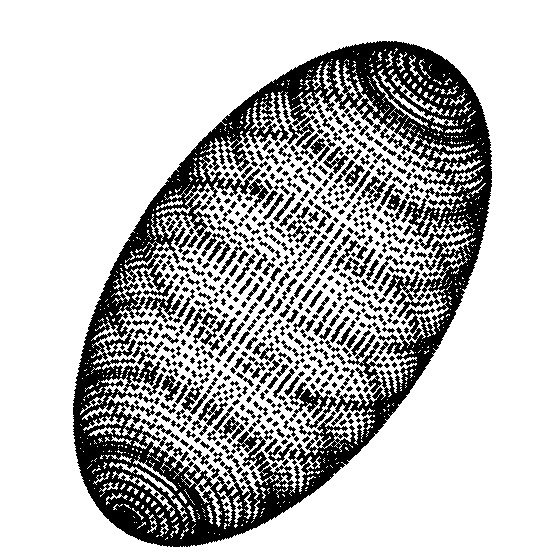}
&\mbox{}\hspace{10mm}\mbox{}&
\includegraphics[width=32mm]{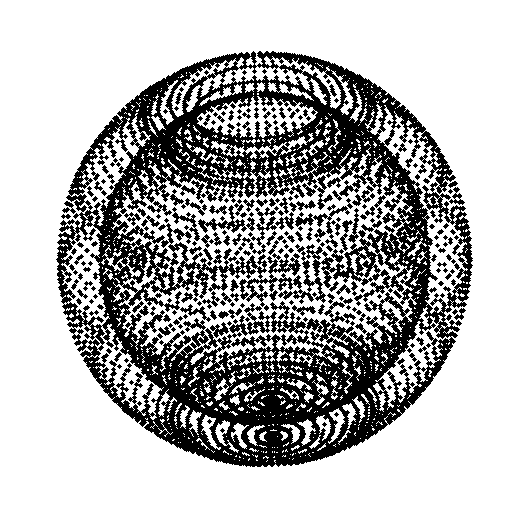}
&\mbox{}\hspace{10mm}\mbox{}&
\includegraphics[width=32mm]{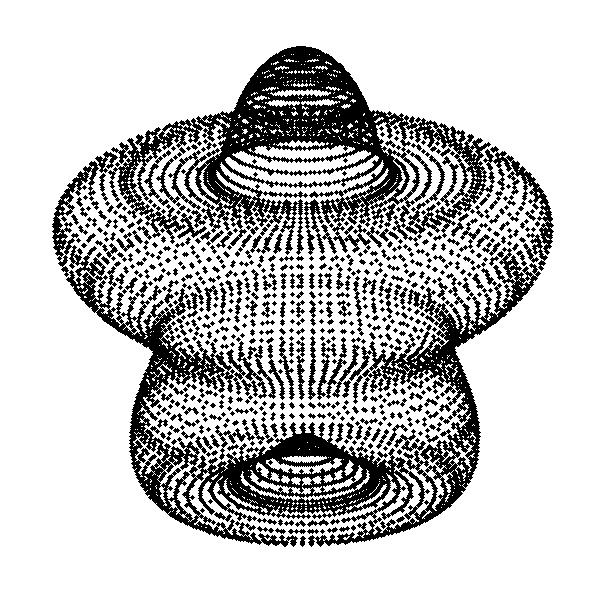}
\\
\includegraphics[width=32mm]{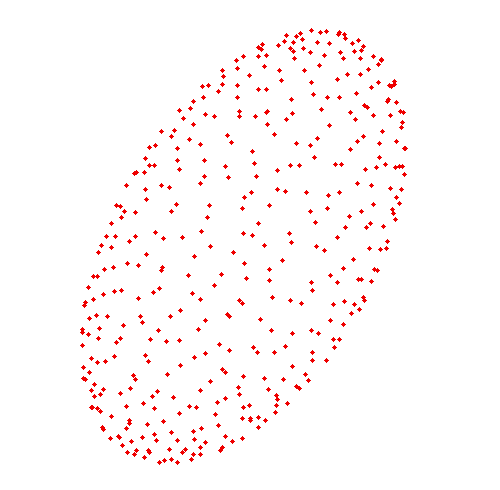}
&\mbox{}\hspace{10mm}\mbox{}&
\includegraphics[width=32mm]{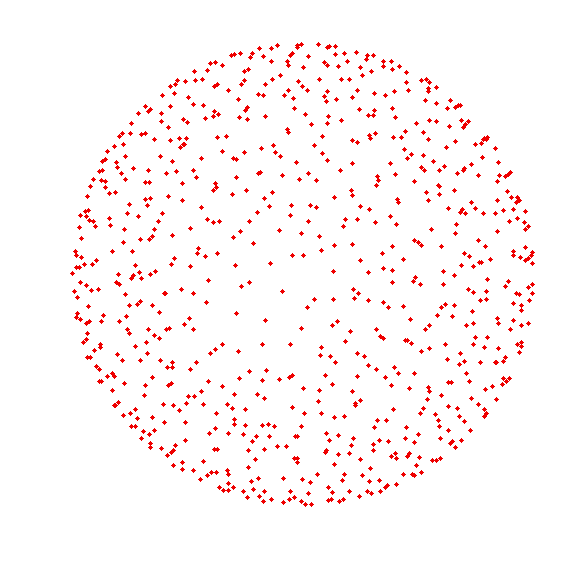}
&\mbox{}\hspace{10mm}\mbox{}&
\includegraphics[width=36mm]{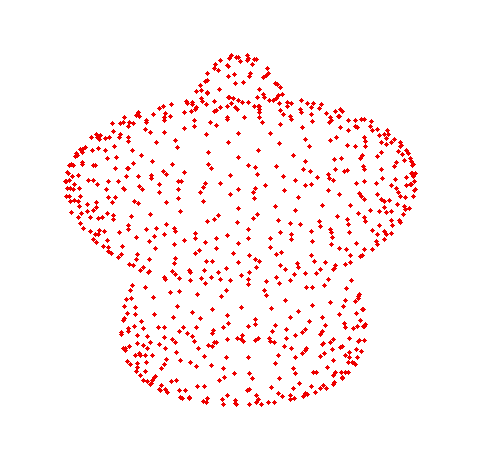}
\\
(a) && (b) && (c)
\end{tabular}
\caption{Example of skeletonization of three different scatterers before and after compression. With $10\,100$ original discretization points (denoted by black dots), after compression (a) for an ellipsoid, only $435$ points survive (denoted by red dots); (b) for a bowl-shaped cavity domain, only $826$ points survive; (c) for a starfish-shaped cavity,  only $803$ points survive.}
\label{fig:compress_example}
\end{figure}


\begin{table}[htbp]
\centering
\begin{tabular}{|r|r|r|c|c|c|c|c|}
\hline
$N$ &    $n$  &  $N_{\textrm{compressed}}$ & $k$  & $T_{\textrm{compress}}$ & $I$ & $T_{\textrm{solve}}$ & $E^{\textrm{rel}}_{\infty}$\\
\hline
127\,500 & $50 \times 51$   & 30\,286 & (411,797,746)  & 3.33e+01 & 18 & 3.85e+01 & 3.042e-05\\
255\,000 & $100\times 51$  & 33\,876  & (434,824,805) & 7.00e+01& 19 & 4.25e+01 & 1.458e-05\\
510\,000 & $200 \times 51$ & 35\,042  & (449,847,838) & 1.46e+02 & 19 & 4.26e+01 & 1.285e-05\\
\hline
252\,500   & $50 \times 101$   & 32\,186 & (413,795,752) & 6.66e+01& 19 & 3.94e+01 & 3.008e-05\\
505\,000   & $100 \times 101$ & 33\,894  & (435,826,803) & 1.40e+02 & 19 & 4.04e+01 & 9.134e-06 \\
1\,010\,000 & $200 \times 101$ & 35\,094 & (451,846,840) & 3.20e+02   & 19 & 4.12e+01 & 5.287e-07\\
\hline
502\,500   & $50 \times 201$   & 32\,286  & (414,797,754) & 1.33e+02 & 19 & 3.98e+01 & 3.013e-05\\
1\,050\,000   & $100 \times 201$ & 33\,798 & (437,830,802)  & 3.00e+02 & 19 & 4.06e+01 & 9.130e-06 \\
2\,010\,000 & $200 \times 201$ & 35\,194   & (453,848,842) & 5.78e+02 & 19 & 4.21e+01 & 4.725e-08\\
\hline
\end{tabular}
\vspace{2mm}
\caption{Example 5: exterior Laplace problem solved on the domain in Figure \ref{fig:diff_scatter} using accelerated scheme. }
\label{tab:lap_comp}
\end{table}

\begin{table}[htbp]
\centering
\begin{tabular}{|r|r|c|c|c|c|}
\hline
$N$ &    $n$  & $T_{\textrm{pre}}$ &  $I$  & $T_{\textrm{solve}}$ &  $E^{\textrm{rel}}_{\infty}$\\
\hline
127\,500 & $50 \times 51$ & 2.29e+00   &  18 & 1.52e+02  &2.908e-05\\
255\,000 & $100\times 51$  & 4.70e+00 & 18 &  2.94e+02   &2.329e-05\\
510\,000 & $200 \times 51$  & 1.22e+01 & 18 & 5.85e+02  & 2.034e-05\\
\hline
252\,500   & $50 \times 101$ & 3.23e+00   & 19 & 2.85e+02 &3.677e-05\\
505\,000   & $100 \times 101$  & 7.08e+00 & 19 & 5.29e+02&1.705e-06 \\
1\,010\,000 & $200 \times 101$ & 1.93e+01 & 19 & 1.06e+03  &4.128e-07\\
\hline
502\,500   & $50 \times 201$    & 5.07e+00  & 19 & 5.02e+02& 3.674e-05\\
1\,050\,000   & $100 \times 201$ & 1.28e+01 & 19 & 9.88e+02  & 1.673e-06 \\
2\,010\,000 & $200 \times 201$ & 3.63e+01 & 19 & 2.07e+03  &1.568e-08\\
\hline
\end{tabular}
\vspace{2mm}
\caption{Example 5: exterior Laplace problem solved on the domain in Figure \ref{fig:diff_scatter} without using accelerated scheme.}
\label{tab:lap_nocomp}
\end{table}

\subsubsection{Example 6}
\begin{figure}[htbp]
\centering
\includegraphics[width=0.45\textwidth]{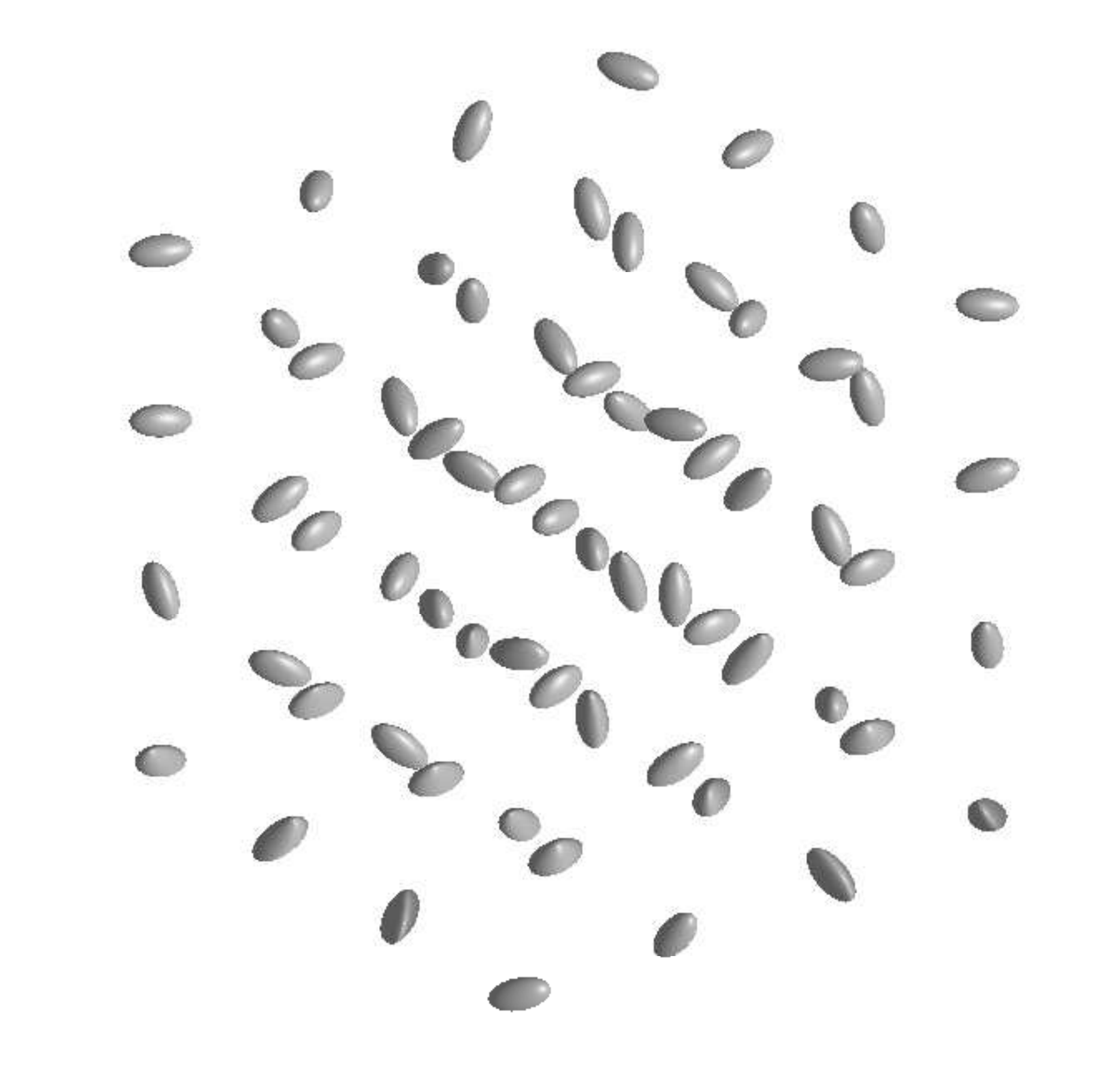}
\caption{Domain contains 64 randomly oriented ellipsoids, where the minimal distance between any two is $6.0$.}
\label{fig:helm_64ellp}
\end{figure}

The accelerated scheme is applied to solve Helmholtz equation on domain containing 64 randomly
placed ellipsoids depicted in Figure \ref{fig:helm_64ellp}. The minimal distance between any
two bodies is $6.0$. Each ellipsoid is 5 wavelength in diameter. We compare the results with and
without accelerated scheme in Table \ref{tab:helm_comp_ellp} and \ref{tab:helm_nocomp_ellp}.
Again, we see that the acceleration scheme leads to very substantial improvements in the solve time.

\begin{table}[htbp]
\centering
\begin{tabular}{|r|r|r|c|c|c|c|c|}
\hline
$N$ &    $n$  &  $N_{\textrm{compressed}}$  &  $k$ & $T_{\textrm{compressed}}$ & $I$ & $T_{\textrm{solve}}$ & $E^{\textrm{rel}}_{\infty}$\\
\hline
80\,000 & $50\times 25$     & 61\,184 & 956   & 1.92e+01 & 28 & 4.42e+03 & 2.339e-02 \\
160\,000& $ 100\times 25$ & 75\,648 & 1182 &  6.58e+01 & 29 & 4.79e+03 & 8.656e-03\\
\hline
163\,200 & $50 \times 51$   & 87\,744   & 1371 & 8.50e+02   & 29 & 4.92e+03 & 2.798e-04\\
326\,400 & $100\times 51$  & 100\,288 & 1567 & 2.83e+03 & 30 & 5.25e+03 & 5.892e-05\\
652\,800 & $200 \times 51$ & 105\,216  & 1644 & 9.06e+02 & 30 & 5.51e+03 & 6.056e-05\\
\hline
323\,200   & $50 \times 101$   & 91\,648   & 1432 & 2.40e+02 & 30 & 5.09e+03 & 9.485e-06\\
646\,400   & $100 \times 101$ & 102\,400 & 1552 & 8.55e+02 & 31 & 5.50e+03 & 2.150e-07 \\
1\,292\,800 & $200 \times 101$ & 106\,944 & 1671 & 2.91e+03 & 31 & 5.73e+03 & 8.441e-08\\
\hline
\end{tabular}
\vspace{2mm}
\caption{Example 6: exterior Helmholtz problem solved on the domain in Figure \ref{fig:helm_64ellp} using accelerated scheme. Each ellipsoid is 5 wavelength in diameter.}
\label{tab:helm_comp_ellp}
\end{table}

\begin{table}[htbp]
\centering
\begin{tabular}{|r|r|c|c|c|c|}
\hline
$N$ &    $n$  & $T_{\textrm{init}}$ &  $I$  & $T_{\textrm{solve}}$ &  $E^{\textrm{rel}}_{\infty}$\\
\hline
80\,000  & $50\times 25$     & 4.41e-01  & 28  & 3.60e+03 & 7.009e-03 \\
160\,000& $ 100\times 25$ & 8.44e-02 &  28  & 5.69e+03 & 5.755e-03\\
\hline
163\,200 & $50 \times 51$   & 8.22e-01   & 28 & 5.78e+03 & 1.239e-04\\
326\,400 & $100\times 51$  & 1.65e+00 & 29  & 8.75e+03 & 4.806e-05\\
652\,800 & $200 \times 51$ & 3.36e+00 & 29  & 1.54e+04 & 5.552e-05\\
\hline
323\,200   & $50 \times 101$   & 1.58e+00 & 29 & 8.64e+03 & 8.223e-06\\
646\,400   & $100 \times 101$ & 3.24e+00 & 29 & 1.69e+04  & 1.354e-07 \\
1\,292\,800 & $200 \times 101$ & 6.67e+00 & 29 & 3.01e+04  &  2.823e-08\\
\hline
\end{tabular}
\vspace{2mm}
\caption{Example 6: exterior Helmholtz problem solved on the domain in Figure \ref{fig:helm_64ellp} without using accelerated scheme. Each ellipsoid is 5 wavelength in diameter.}
\label{tab:helm_nocomp_ellp}
\end{table}

\subsubsection{Example 7}
The accelerated scheme is applied to solve Helmholtz equation on domain in Figure \ref{fig:diff_scatter}.
Each scatterer is two wavelength in diameter. Results are presented in Table \ref{tab:helm_comp} and
Table \ref{tab:helm_nocomp}, and confirm our claim regarding the efficiency of the acceleration scheme.
Note that in Table \ref{tab:helm_nocomp}, due to limitation of the memory, only estimations of the run time are reported when four million discretization nodes were used.

\begin{table}[htbp]
\centering
\begin{tabular}{|r|r|r|c|c|c|c|c|}
\hline
$N$ &    $n$  &  $N_{\textrm{compressed}}$  &  $k$ & $T_{\textrm{compressed}}$ & $I$ & $T_{\textrm{solve}}$ & $E^{\textrm{rel}}_{\infty}$\\
\hline
252\,500 & $50\times 101$  & 53\,390 & (775,1254,1211) & 2.26e+02 & 52 & 2.48e+03 & 4.941e-03 \\
505\,000 & $100\times 101$ & 57\,934 & (823.1358,1337) & 5.17e+02 & 53 & 2.72e+03 & 2.026e-03\\
1\,010\,000& $200\times 101$ & 60\,512 & (856,1420,1399) & 1.14e+03 & 54 & 2.89e+03 & 4.865e-04\\
\hline
502\,500   & $50 \times 201$  & 54\,538   & (789,1283,1238) & 4.89e+02 & 53 & 2.63e+03 & 9.276e-03\\
1\,005\,000  & $100 \times 201$ & 59\,036   & (838,1384,1363) & 1.10e+02 & 54 & 2.90e+03 & 4.392e-03 \\
2\,010\,000  & $200 \times 201$ & 61\,488 & (872,1443,1419) & 2.70e+03 & 56 & 3.10e+03 & 7.709e-06\\
4\,020\,000 & $400 \times 201$ & 61\,664 & (888,1428,1427) & 1.50e+04 & 57 & 3.31e+03& 1.856e-06 \\
\hline
2\,005\,000 & $100 \times 401$ & 60\,106 & (853,1409,1388) &  2.58e+03 & 56 & 3.04e+03 & 9.632e-04\\
4\,010\,000 & $200 \times 401$ & 61\,818 & (885,1441,1427) &1.54e+04  & 57& 3.32e+03& 2.452e-07\\
\hline
\end{tabular}
\vspace{2mm}
\caption{Example 7: exterior Helmholtz problem solved on the domain in Figure \ref{fig:diff_scatter} using accelerated scheme. Each scatterer is 2 wavelength in diameter.}
\label{tab:helm_comp}
\end{table}

\begin{table}[htbp]
\centering
\begin{tabular}{|r|r|c|c|c|c|}
\hline
$N$ &    $n$  & $T_{\textrm{init}}$ &  $I$  & $T_{\textrm{solve}}$ &  $E^{\textrm{rel}}_{\infty}$\\
\hline
252\,500 & $50\times 101$   & 5.33e+00  & 50  & 1.01e+04 & 3.211e-03 \\
505\,000 & $100\times 101$  & 1.07e+01  & 50  & 2.04e+04 & 2.260e-03\\
1\,010\,000& $200\times 101$  & 2.21e+01  & 51  & 4.16e+04 & 8.211e-04\\
\hline
502\,500   & $50 \times 201$ &   1.02e+01 & 51  & 2.15e+04 & 8.273e-03   \\
1\,005\,000  & $100 \times 201$   & 2.01e+01& 51 & 4.20e+04 & 3.914e-03\\
2\,010\,000  & $200 \times 201$  & 3.90e+01  & 51  & 8.42e+04 & 5.044e-06 \\
4\,020\,000 & $400 \times 201$ & -- & -- & $\sim$ 48h& -- \\
\hline
2\,005\,000 & $100\times 401$ & 3.89e+01 & 51 & 8.30e+04 & 4.244e-04\\
4\,010\,000 & $200 \times 401$& -- & -- & $\sim$ 48h & --\\
\hline
\end{tabular}
\vspace{2mm}
\caption{Example 7: exterior Helmholtz problem solved on the domain in Figure \ref{fig:diff_scatter} without using accelerated scheme. Each scatterer is 2 wavelength in diameter.}
\label{tab:helm_nocomp}
\end{table}

\section{Conclusions and Future work}
\label{sec:conclu}

We have presented a highly accurate numerical scheme for solving acoustic
scattering problems on domains involving multiple scatterers in three dimensions,
under the assumption that each scatterer is axisymmetric. The algorithm relies
on a boundary integral equation formulation of the scattering problem, combined
with a highly accurate Nystr\"om discretization technique. For each scatterer, a
scattering matrix is constructed via an explicit inversion scheme. Then these
individual scattering matrices are used as a block-diagonal pre-conditioner to
GMRES to solve the very large system of linear equations. The Fast Multiple Method
is used to accelerate the evaluation of all inter-body interactions. Numerical
experiments show that while the block-diagonal pre-conditioner does not make almost
any different for ``zero-frequency'' scattering problems (governed by Laplace's
equation), it dramatically improves the convergence speed at intermediate frequencies.

Furthermore, for problems where the scatterers are well-separated, we present
an accelerated scheme capable of solving even very large scale problems to high
accuracy on a basic personal work station. In one numerical
examples in Section \ref{sec:num}, the numbers of degrees of freedom required to
solve the Laplace equation to eight digits of accuracy on a complex geometry could
be reduced by a factor of $57$ resulting in a reduction of the total computation
time from 35 minutes to $10$ minutes (9 minutes for compression and 42 seconds for
solving the linear system).
For a Helmholtz problem the reduction of computation time is even more significant:
the numbers of degrees of freedom to reach seven digits
of accuracy was in one example reduced by a factor of $65$; consequently the overall computation
time is reduced from $48$ hours to $5$ hours (4 hours for compression and 1 hour for solving the linear system).

The scheme presented assumes that each scatterer is rotationally symmetric; this
property is used both to achieve higher accuracy in the discretization, and to
accelerate all computations (by using the FFT in the azimuthal direction).
It appears conceptually straight-forward to use the techniques of
\cite{2012_bremer_nystrom,2008_helsing_corner_BIE} to generalize the method
presented to handle scatterers with edges (generated by ``corners'' in the
generating curve). The idea is to use local refinement to resolve the singular
behavior of solutions near the corner, and then eliminate the added ``superfluous''
degrees of freedom added by the refinement via a local compression technique,
see \cite{gillman2014short}.

\noindent\\
\textbf{\textit{Acknowledgements:}}
The research reported was supported by the National Science Foundation under
contracts 1320652, 0941476, and 0748488, and by the Defense Advanced Projects
Research Agency under the contract N66001-13-1-4050.

\bibliography{main_bib,refs}
\bibliographystyle{amsplain}


\end{document}